\documentclass[letterpaper, 10 pt, conference]{ieeeconf}  

	\IEEEoverridecommandlockouts                              
	\PassOptionsToPackage{noend}{algpseudocode}
	\let\labelindent\relax

	\usepackage{amsthm}
	\usepackage{myPreamble}
	
	\renewcommand\eg{\textit{e.g.}\ }
	\DeclareMathOperator\sym{Sym}
	\newcommand\LL{\mathcal L}
	\newcommand\M{\mathcal M}
	\newcommand\x{\text{\textsf{\textsc x}}}
	\newcommand\y{\text{\textsf{\textsc y}}}
	
	\makeatletter
		\def\Sx{%
			\Sigma%
			\@ifnextchar_{\Sx@}{_\x}%
		}
		\def\Sx@_#1{%
			_{\x,#1}%
		}
		\def\Sy{%
			\Sigma%
			\@ifnextchar_{\Sy@}{_\y}%
		}
		\def\Sy@_#1{%
			_{\y,#1}%
		}
	\makeatother

		\makeatletter
			\renewcommand\alglinenumber[1]{\footnotesize\textbf{\thealgorithm}.\oldstylenums{\arabic{ALG@line}}:}
			\renewcommand\theALG@line{\thealgorithm.\oldstylenums{\arabic{ALG@line}}}
			\providecommand\theHALG@line{\thealgorithm.\arabic{ALG@line}}

			\let\OLDIf\If
			\let\OLDEndIf\EndIf
			\let\OLDState\State
			\newcommand\@@If[1]{{\keywordfont if\ }#1\ {\keywordfont then}\def\State{\ }}
			\newcommand\@@EndIf{\ {\keywordfont end if}\def\State{\OLDState}\ignorespaces}
			\renewcommand\If{\@ifstar\@@If\OLDIf}
			\renewcommand\EndIf{\@ifstar\@@EndIf\OLDEndIf}
		\makeatother
		
		\newcommand\keywordfont{\scshape}
		\algrenewcommand\algorithmicend{{\keywordfont end}}
		\algrenewcommand\algorithmicdo{{\keywordfont do}}
		\algrenewcommand\algorithmicwhile{{\keywordfont while}}
		\algrenewcommand\algorithmicfor{{\keywordfont for}}
		\algrenewcommand\algorithmicforall{{\keywordfont for all}}
		\algrenewcommand\algorithmicloop{{\keywordfont loop}}
		\algrenewcommand\algorithmicrepeat{{\keywordfont repeat}}
		\algrenewcommand\algorithmicuntil{{\keywordfont until}}
		\algrenewcommand\algorithmicprocedure{{\keywordfont procedure}}
		\algrenewcommand\algorithmicfunction{{\keywordfont function}}
		\algrenewcommand\algorithmicif{{\keywordfont if}}
		\algrenewcommand\algorithmicthen{{\keywordfont then}}
		\algrenewcommand\algorithmicelse{{\keywordfont else}}
		\algrenewcommand\algorithmicrequire{{\keywordfont Require:}}
		\algrenewcommand\algorithmicensure{{\keywordfont Ensure:}}
		\algrenewcommand\algorithmicreturn{{\keywordfont return}}
	
		\renewcommand\algorithmicrequire{\fillwidthof[l]{{\keywordfont Initialize}}{{\keywordfont Require}}}
		\renewcommand\algorithmicensure{\fillwidthof[l]{{\keywordfont Initialize}}{{\keywordfont Ensure}}}
		\renewcommand\algorithmicprovide{\fillwidthof[l]{{\keywordfont Initialize}}{{\keywordfont Provide}}}
		\renewcommand\algorithmicinitialize{\fillwidthof[l]{{\keywordfont Initialize}}{{\keywordfont Initialize}}}

\title{\LARGE \bf QPALM: A Newton-type Proximal Augmented Lagrangian Method\texorpdfstring{\\}{ }for Quadratic Programs}
\author{%
	Ben Hermans,\texorpdfstring{\textsuperscript1}{}
	Andreas Themelis\texorpdfstring{\textsuperscript2}{} and
	Panagiotis Patrinos\texorpdfstring{\textsuperscript2}{}%
	\texorpdfstring{\thanks{%
		\textsuperscript1Ben Hermans is with the Department of Mechanical Engineering, KU Leuven, and DMMS lab, Flanders Make, Leuven, Belgium.
		His research benefits from KU Leuven-BOF PFV/10/002 Centre of Excellence: Optimization in Engineering (OPTEC), from project G0C4515N of the Research Foundation - Flanders (FWO - Flanders), from Flanders Make ICON: Avoidance of collisions and obstacles in narrow lanes, and from the KU Leuven Research project C14/15/067: B-spline based certificates of positivity with applications in engineering.%
		\newline
		\textsuperscript2Andreas Themelis and Panagiotis Patrinos are with the \TheAddressKU.
		This work was supported by the \emph{Research Foundation Flanders (FWO)} research projects G086518N and G086318N;
		\emph{Research Council KU Leuven} C1 project No. C14/18/068;
		\emph{Fonds de la Recherche Scientifique --- FNRS and the Fonds Wetenschappelijk Onderzoek --- Vlaanderen} under EOS project no 30468160 (SeLMA).%
		\newline
		{\tt
			\{\href{mailto:ben.hermans2@kuleuven.be}{ben.hermans2},%
			\href{mailto:andreas.themelis@kuleuven.be}{andreas.themelis},%
			\href{mailto:panos.patrinos@kuleuven.be}{panos.patrinos}\}%
			\href{mailto:ben.hermans2@kuleuven.be,andreas.themelis@kuleuven.be,panos.patrinos@kuleuven.be}{@kuleuven.be}%
		}%
	}%
}{}}

\begin{document}

	\maketitle
	\thispagestyle{empty}
	\pagestyle{empty}
	\begin{abstract}
		We present a proximal augmented Lagrangian based solver for general convex quadratic programs (QPs), relying on semismooth Newton iterations with exact line search to solve the inner subproblems.
		The exact line search reduces in this case to finding the zero of a one-dimensional monotone, piecewise affine function and can be carried out very efficiently.
		Our algorithm requires the solution of a linear system at every iteration, but as the matrix to be factorized depends on the active constraints, efficient sparse factorization updates can be employed like in active-set methods.  
		Both primal and dual residuals can be enforced down to strict tolerances and otherwise infeasibility can be detected from intermediate iterates.
		A C implementation of the proposed algorithm is tested and benchmarked against other state-of-the-art QP solvers for a large variety of problem data and shown to compare favorably against these solvers.
	\end{abstract}

\section{Introduction}
	This paper deals with convex quadratic programs (QPs)
	\[\tag{QP}\label{eq:QP}
		\minimize_{x\in\R^n}\tfrac12\innprod{x}{Qx}+\innprod qx
	\quad\stt{}
		\ell\leq Ax\leq u,
	\]
	where \(Q\in\R^{n\times n}\) is symmetric and positive semidefinite, \(q\in\R^n\), \(A\in\R^{m\times n}\), and \(\ell,u\in\R^m\).
	
	Efficiently and reliably solving QPs is a key challenge in optimization \cite[\S16]{nocedal2006numerical}. QPs cover a wide variety of applications and problem classes, such as portfolio optimization, support vector machines, sparse regressor selection, real-time linear model predictive control (MPC), etc.
	QPs also often arise as subproblems in general nonlinear optimization techniques such as sequential quadratic programming.
	Therefore, substantial research has been performed to develop robust and efficient QP solvers.
	State-of-the-art solvers are typically based on interior point methods, such as MOSEK \cite{mosek} and Gurobi \cite{gurobi2018gurobi}, or active-set methods, such as qpOASES \cite{ferreau2008online}.
	Both methods have their advantages and disadvantages.
	Interior point methods typically require few but expensive iterations, involving the solution of a linear system at every iteration.
	In contrast, active-set methods require more but cheaper iterations, as the linear system changes only slightly and low-rank factorization updates can be used instead of having to factorize from scratch.
	As a result, active-set methods can be efficiently warm started when solving a series of similar QPs, which is common in applications such as MPC, whereas interior-point methods in general do not have this capability.
	Recently, proximal algorithms, also known as operator splitting methods \cite{parikh2014proximal}, have experienced a resurgence in popularity.
	Relying only on first-order information, such methods have as their advantage operational simplicity and cheap iterations, but they may exhibit slow asymptotic convergence for poorly scaled problems.
	The recently proposed OSQP solver \cite{stellato2017osqp}, based on the alternating direction method of multipliers (ADMM), addresses this crucial issue by means of a tailored offline scaling that performs very well on some QP problems, although a more thorough benchmarking confirms the known limitations of proximal algorithms.
	Indeed, parameter tuning is typically set before execution, with possibly minor online adjustments, and as such operator splitting methods do not take into account curvature information about the problem that can greatly speed up convergence.
	
	In this work we propose a novel and reliable QP solver based on the augmented Lagrangian method (ALM) \cite{bertsekas1999constrained} and in particular proximal ALM \cite{rockafellar1976augmented}, which is efficient and robust against ill conditioning.
	ALM involves penalizing the constraints and solving a series of unconstrained minimizations, where fast smooth optimization techniques can be employed.
	QPs turn out to be particularly amenable for this approach, as optimal stepsizes of exact line searches are available in closed form, similiar to what was shown in previous work \cite{patrinos2011global}, resulting in an extremely fast minimization strategy.
	In each unconstrained optimization, the iterates rely on the solution of a linear system, dependent on the set of active constraints.
	As such, similarly to active-set methods, our iterates can benefit from low-rank factorization update techniques and are therefore cheaper than in interior point methods.
	However, in contrast to active-set methods and more in the flavor of other algorithms such as the dual gradient projection method of \cite{axehill2008dual}, our method allows for substantial changes in the active set at each iteration and converges therefore much faster on average.
	To some extent, our algorithm strikes a balance between interior point and active-set methods. In regard to state-of-the-art ALM solvers, not much research is available. The authors of \cite{gilbert2014oqla} presented an ALM solver for QPs, OQLA/QPALM, which solves the inner minimization problems using a combination of active-set and gradient projection methods. 
	Although the authors of \cite{Dhingra2017second} discuss second order updates based on a generalized Hessian, similarly to our method their approach deals with convex (composite) optimization problems in general, and as such it can not make use of the efficient factorization updates and optimal stepsizes that are key to the efficiency of our algorithm.
	Finally, a link with proximal algorithms can also be established, in the sense that the form of an iterate of unconstrained minimization is very similar to what is obtained by applying ADMM to \eqref{eq:QP}. Differently from an ADMM approach where parameters are set before the execution, in our method penalties are adjusted online. Moreover, penalty matrices are used instead of scalars to allow for more flexibility in penalizing the constraints. As a result of the online adjustments, our algorithm requires in general fewer, although slightly more expensive, iterations.
	
	The contribution of this paper is QPALM, a full-fledged convex QP solver based on proximal ALM.
	We describe the relevant theory and detail the algorithmic steps for both the outer and inner minimization procedures.
	We also provide an open-source C implementation of the algorithm,\footnote{\url{https://github.com/Benny44/QPALM}} which we benchmark with state-of-the-art QP solvers, showing that QPALM can compete with and very often outperform them both in runtime and robustness in regard to problem scaling.
	
	The remainder of this paper is structured as follows.
	\Cref{sec:Preliminaries} introduces notation used in the paper and underlying theoretical concepts.
	\Cref{sec:PALM-QP} outlines the application of the proximal ALM to \eqref{eq:QP}.
	\Cref{sec:QPALM} discusses in detail the algorithmic steps, regarding both inner and outer minimization, of QPALM.
	\Cref{sec:Implementation} addresses some of the crucial implementation details that contribute to making QPALM an efficient and robust solver.
	\Cref{sec:Numerical simulations} presents simulation results on QPs of varying sizes and problem conditioning, benchmarking QPALM's performance against state-of-the-art solvers.
	Finally, \Cref{sec:Conclusion} draws concluding remarks and mentions future directions of research.


	\section{Preliminaries}
		We now introduce some notational conventions and briefly list some known facts needed in the paper.
		The interested reader is referred to \cite{bauschke2017convex} for an extensive discussion.
		\label{sec:Preliminaries}
		\subsection{Notation and known facts}
	We denote the extended real line by \(\Rinf\coloneqq\R\cup\set{\infty}\).
	The scalar product on \(\R^n\) is denoted by \(\innprod{{}\cdot{}}{{}\cdot{}}\).
	With $[x]_+\coloneqq\max\set{x,0}$ we indicate the positive part of vector \(x\in\R^n\), meant in a componentwise sense.
	A sequence of vectors \(\seq{x^k}\) is said to be summable if \(\sum_{k\in\N}\|x^k\|<\infty\).
	
	With \(\sym(\R^n)\) we indicate the set of symmetric \(\R^{n\times n}\) matrices, and with \(\sym_+(\R^n)\) and \(\sym_{++}(\R^n)\) the subsets of those which are positive semidefinite and positive definite, respectively.
	Given \(\Sigma\in\sym_{++}(\R^n)\) we indicate with
	\(
		\|{}\cdot{}\|_\Sigma
	\)
	the norm on \(\R^n\) induced by \(\Sigma\), namely
	\(
		\|x\|_\Sigma
	{}\coloneqq{}
		\sqrt{\innprod{x}{\Sigma x}}
	\).
	
	Given a nonempty closed convex set \(C\subseteq\R^n\), with \(\proj_C(x)\) we indicate the projection of a point \(x\in\R^n\) onto \(C\), namely \(\proj_C(x)=\argmin_{y\in C}\|y-x\|\) or, equivalently, the unique point \(z\in C\) satisfying the inclusion
	\begin{equation}\label{eq:proj}
		x-z\in\ncone_C(z),
	\end{equation}
	where
	\(
		\ncone_C(z)
	{}\coloneqq{}
		\set{v\in\R^n}[\innprod{v}{z-z'}\leq 0~\forall z'\in C]
	\)
	is the normal cone of \(C\) at \(z\).
	\(\dist(x,C)\) and \(\dist_\Sigma(x,C)\) denote the distance from \(x\) to set \(C\) in the Euclidean norm and in that induced by \(\Sigma\), respectively, while \(\indicator_C\) is the indicator function of set \(C\), namely
	\(
		\indicator_C(x)=0
	\)
	if \(x\in C\) and \(\infty\) otherwise.
		\subsection{Convex functions and monotone operators}
	The \DEF{Fenchel conjugate} of a proper closed convex function \(\func\varphi{\R^n}{\Rinf}\) is the convex function
	\(
		\func{\conj\varphi}{\R^n}{\Rinf}
	\)
	defined as
	\(
		\conj\varphi(y)
	{}={}
		\sup_x{\innprod xy-\varphi(x)}
	\).
	The subdifferential of \(\varphi\) at \(x\in\R^n\) is
	\(
		\partial\varphi(x)
	{}\coloneqq{}
		\set{v\in\R^n}[
			\varphi(x')\geq\varphi(x)+\innprod{v}{x'-x}
			~
			\forall x'\in\R^n
		]
	\).
	Having \(y\in\partial\varphi(x)\) is equivalent to \(x\in\partial\conj\varphi(y)\).
	
	A point-to-set mapping \(\ffunc{\M}{\R^n}{\R^n}\) is \DEF{monotone} if
	\(
		\innprod{x-x'}{\xi-\xi'}
	{}\geq{}
		0
	\)
	for all \(x,x'\in\R^n\), \(\xi\in\M(x)\) and \(\xi'\in\M(x')\).
	It is \DEF{maximally monotone} if, additionally, there exists no monotone operator \(\M'\neq\M\) such that \(\M(x)\subseteq\M'(x)\) for all \(x\in\R^n\).
	The \DEF{resolvent} of a maximally monotone operator \(\M\) is the single-valued (in fact, Lipschitz-continuous) mapping \((\id+\M)^{-1}\), where
	\(
		(\id+\M)^{-1}(x)
	\)
	is the unique point \(\bar x\in\R^n\) such that \(x-\bar x\in\M(\bar x)\).
	\(\zer\M\coloneqq\set{x}[0\in\M(x)]\) denotes the zero-set of \(\M\), and for a linear mapping \(\Sigma\), \(\Sigma\M\) is the operator defined as \(\Sigma M(x)\coloneqq\set{\Sigma y}[y\in\M(x)]\).
	
	The subdifferential of a proper convex lower semicontinuous function \(\varphi\) is maximally monotone, and its resolvent is the proximal mapping
	\(
		\prox_\varphi\coloneqq(\id+\partial\varphi)^{-1}
	\).

	\section{Proximal ALM for QPs} 
		\label{sec:PALM-QP}
	This section is devoted to establishing the theoretical ground in support of the proposed \Cref{alg:QPALM}.
	We will show that our scheme amounts to proximal ALM and derive its convergence guarantees following the original analysis in \cite{rockafellar1976augmented}, here generalized to account for scaling matrices (as opposed to scalars) and by including the possibility of having different Lagrangian and proximal weights.
	We start by observing that problem \eqref{eq:QP} can equivalently be expressed as
	\begin{equation}\label{eq:P}
		\minimize_{x\in\R^n}f(x)+g(Ax),
	\end{equation}
	where
	\(
		f(x)
	{}\coloneqq{}
		\tfrac12\innprod{x}{Qx}+\innprod qx,
	\)
	\(
		g(z)
	{}\coloneqq{}
		\indicator_C(z)
	\)
	and
	\(
		C
	{}\coloneqq{}
		\set{z\in\R^m}[\ell\leq z\leq u]
	\).
	The KKT conditions of \eqref{eq:P} are
	\[
		0\in\M(x,y)
	{}\coloneqq{}
		\binom{
			\nabla f(x)+\trans Ay
		}{
			-Ax+\partial\conj g(y)
		}.
	\]
	Let \(\mathcal V_\star\coloneqq\zer\M\) be the set of primal-dual solutions.
	Since \(\M\) is maximally monotone, as first observed in \cite{rockafellar1976augmented} one can find KKT-optimal primal-dual pairs by recursively applying the resolvent of \(c_k\M\), where \(\seq{c_k}\) is an increasing sequence of strictly positive scalars.
	This scheme is known as proximal point algorithm (PPA) \cite{rockafellar1976augmented}.
	We now show that these scalars can in fact be replaced by positive definite matrices.
	\begin{thm}\label{thm:SigmaM}%
		Suppose that \eqref{eq:P} has a solution.
		Starting from \((x^0,y^0)\in\R^n\times\R^m\), let \(\seq{x^k,y^k}\) be recursively defined as
		\begin{equation}\label{eq:SigmaM}
			(x^{k+1},y^{k+1})=(\id+\Sigma_k\M)^{-1}(x^k,y^k)+\varepsilon^k
		\end{equation}
		for a summable sequence \(\seq{\varepsilon^k}\), and where \(\Sigma_k\coloneqq\binom{\Sx_k~\hphantom{\Sy_k}}{\hphantom{\Sx_k}~\Sy_k}\) for some \(\Sx_k\in\sym_{++}(\R^n)\) and \(\Sy_k\in\sym_{++}(\R^m)\).
		If \(\Sigma_k\preceq\Sigma_{k+1}\preceq\Sigma_\infty\in\sym_{++}(\R^n\times\R^m)\) holds for all \(k\), then \(\seq{x^k,y^k}\) converges to a KKT-optimal pair for \eqref{eq:P}.
		\begin{proof}
			We start by observing that for all \(k\) it holds that \(\zer(\Sigma_k\M)=\zer(\M)=\mathcal V_\star\) and that \(\Sigma_k\M\) is maximally monotone with respect to the scalar product induced by \(\Sigma_k^{-1}\).
			The resolvent \((\id+\Sigma_k\M)^{-1}\) is thus firmly nonexpansive in that metric (see \cite[Prop. 23.8 and Def. 4.1]{bauschke2017convex}):
			that is, denoting \(v^k\coloneqq(x^k,y^k)\) and \(\tilde v^{k+1}\coloneqq v^{k+1}-\varepsilon^k=(\id+\Sigma_k\M)^{-1}(v^k)\),
			\begin{equation}\label{eq:Fejer}
				\|\tilde v^{k+1}-v_\star\|_{\Sigma_k^{-1}}^2
			{}\leq{}
				\|v^k-v_\star\|_{\Sigma_k^{-1}}^2
				{}-{}
				\|\tilde v^{k+1}-v^k\|_{\Sigma_k^{-1}}^2
			\end{equation}
			holds for every \(v_\star\in\mathcal V_\star\).
			Therefore, since
			\(
				\|v^{k+1}-v_\star\|_{\Sigma_{k+1}^{-1}}
			{}\leq{}
				\|v^{k+1}-v_\star\|_{\Sigma_k^{-1}}
			{}\leq{}
				\|\tilde v^{k+1}-v_\star\|_{\Sigma_k^{-1}}
				{}+{}
				\|\varepsilon^k\|_{\Sigma_k^{-1}}
			{}\leq{}
				\|v^k-v_\star\|_{\Sigma_k^{-1}}
				{}+{}
				\|\varepsilon^k\|_{\Sigma_k^{-1}}
			\)
			(where the first inequality owes to the fact that \(\Sigma_{k+1}^{-1}\preceq\Sigma_k^{-1}\)), it follows from \cite[Thm. 3.3]{combettes2013variable} that the proof reduces to showing that any limit point of \(\seq{v^k}\) belongs to \(\mathcal V_\star\).
			
			From \cite[Prop. 3.2(i)]{combettes2013variable} it follows that the sequence is bounded and that \(v^{k+1}-v^k\to0\) as \(k\to\infty\).
			Suppose that a subsequence \(\seq{v^{k_j}}[j\in\N]\) converges to \(v\); then, so do \(v^{k_j+1}\) and
			\(
				\tilde v^{k_j+1}
			{}={}
				v^{k_j+1}-\varepsilon^{k_j}
			{}={}
				(\id+\Sigma_{k_j}\M)^{-1}v^{k_j}
			\).
			We have
			\[
				\Sigma_{k_j}^{-1}(v^{k_j}-\tilde v^{k_j+1})
			{}\in{}
				\M(\tilde v^{k_j+1});
			\]
			since \(\seq{\Sigma_k^{-1}}\) is upper bounded, the left-hand side converges to \(0\), and from outer semicontinuity of \(\M\) \cite[Ex. 12.8(b)]{rockafellar2011variational} it then follows that \(0\in\M(v)\), proving the claim.
		\end{proof}
	\end{thm}
	\begin{subequations}\label{subeq:xy}%
		Let us consider the iterates \eqref{eq:SigmaM} under the assumptions of \Cref{thm:SigmaM}, and let us further assume that \(\Sy_k\) is diagonal for all \(k\).
		Let \((\tilde x^{k+1},\tilde y^{k+1})\coloneqq (x^{k+1},y^{k+1})-\varepsilon^k=(\id+\Sigma_k\M)^{-1}(x^k,y^k)\).
		Equation \eqref{eq:SigmaM} reads
		\begin{align}
			\label{eq:x+}
			0
		{}={}&
			\Sx_k^{-1}(\tilde x^{k+1}-x^k)+\nabla f(\tilde x^{k+1})+\trans A\tilde y^{k+1}
		\\
			\label{eq:y+}
			0
		{}\in{}&
			\tilde y^{k+1}-y^k+\Sy_k(\partial\conj g(\tilde y^{k+1})-A\tilde x^{k+1}).
		\end{align}
		Let \(\prox_g^\Sigma({}\cdot{}) = \argmin_y(g({}\cdot{}) +  \tfrac12\|\cdot {}-{} y\|_{\Sigma^{-1}}^2)\). 
		Then \eqref{eq:y+} is equivalent to
		\begin{align*}
			\tilde y^{k+1}
		{}={} &
			\prox_{\conj g}^{\Sy_k}(y^k+\Sy_k A\tilde x^{k+1})
		\\
		{}={} &
			y^k+\Sy_k A\tilde x^{k+1}
			{}-{}
			\Sy_k\prox_g^{\Sy_k^{-1}}(A\tilde x^{k+1}+\Sy_k^{-1}y^k)
		\\
		{}={} &
			\Sy_k\left(
				A\tilde x^{k+1}+\Sy_k^{-1}y^k
				{}-{}
				\proj_C(A\tilde x^{k+1}+\Sy_k^{-1}y^k)
			\right),
		\end{align*}
		where the second equality follows from the Moreau decomposition \cite[Thm. 14.3(ii)]{bauschke2017convex}, and the third one from the fact that \(\Sy_k\) is diagonal and set \(C\) is separable, hence that the projection on \(C\) with respect to the Euclidean norm and that induced by \(\Sy_k\) coincide.
	\end{subequations}
	Notice that \(\trans A\tilde y^{k+1}\) is the gradient of \(\frac12\dist_{\Sy_k}^2(A{}\cdot{}+\Sy_k^{-1}y^k, C)\) at \(\tilde x^{k+1}\).
	Using this in \eqref{eq:x+}, by introducing an auxiliary variable \(z^k\) we obtain that an (exact) resolvent step \((\tilde x^{k+1},\tilde y^{k+1})=(\id+\Sigma_k\M)^{-1}(x^k,y^k)\) amounts to
	\begin{equation}\label{eq:PALM}
		\begin{cases}[r @{{}={}} l]
			\tilde x^{k+1} & \argmin_x\varphi_k(x)\\
			\tilde z^{k+1} & Z_k(\tilde x^{k+1})\\
			\tilde y^{k+1} & y^k+\Sy_k\bigl(A\tilde x^{k+1}-\tilde z^{k+1}\bigr),
		\end{cases}
	\end{equation}
	where
	\begin{align*}
		Z_k(x)
	{}\coloneqq{} &
		\argmin_{\smash{z\in C}}{
			\tfrac12\|z-(Ax+\Sy_k^{-1}y^k)\|_{\Sy_k}^2
		}
	{}={}
		\proj_C(Ax+\Sy_k^{-1}y^k)
	\\
	\numberthis\label{eq:z}
	{}={} &
		\mbox{\(\mathtight[0.33]
			Ax+\Sy_k^{-1}y^k
			{}+{}
			[\ell-Ax-\Sy_k^{-1}y^k]_+
			{}-{}
			[Ax+\Sy_k^{-1}y^k-u]_+
		\)}
	\end{align*}
	is a Lipschitz-continuous mapping, and
	\begin{equation}\label{eq:phi}
		\varphi_k(x)
	{}\coloneqq{}
		f(x)+\tfrac12\dist_{\Sy_k}^2(Ax+\Sy_k^{-1}y^k,C)
		{}+{}
		\tfrac12\|x-x^k\|_{\Sx_k^{-1}}^2
	\end{equation}
	is (Lipschitz) differentiable and strongly convex with
	\begin{equation}\label{eq:nabla}
		\nabla\varphi_k(x)
	{}={}
		\nabla f(x)
		{}+{}
		\trans A\bigl(
			y^k+\Sy_k(Ax-Z_k(x))
		\bigr)
		{}+{}
		\Sx_k^{-1}(x-x^k).
	\end{equation}
	\begin{rem}[Connection with the proximal ALM~\cite{rockafellar1976augmented}]\label{thm:alg=PALM}%
		The \((x,z)\)-update in \eqref{eq:PALM} can equivalently be expressed as
		\[
			(x^{k+1},z^{k+1})
		{}={}
			\argmin_{(x,z)\in\R^n\times\R^m}\LL_{\Sy_k}(x,z,y^k)+\tfrac12\|x-x^k\|_{\Sx_k^{-1}}^2,
		\]
		where for \(\Sy\in\sym_{++}(\R^m)\)
		\[
			\LL_{\Sy}(x,z,y)
		{}\coloneqq{}
			f(x)+\indicator_C(z)+\innprod{y}{Ax}+\tfrac12\|Ax-z\|_{\Sy}^2
		\]
		is the \(\Sy\)-augmented Lagrangian associated to
		\begin{equation}\label{eq:Pxz}
			\minimize_{x\in\R^n,z\in\R^m}f(x)+\indicator_C(z)
		\quad\stt{}
			Ax=z,
		\end{equation}
		a formulation equivalent to \eqref{eq:P}.
		In fact, the iterative scheme \eqref{eq:PALM} simply amounts to the proximal ALM applied to \eqref{eq:Pxz}.
	\end{rem}

	\section{The QPALM algorithm}
		\label{sec:QPALM}
	\begin{algorithm}[t]
		\algcaption{Quadratic Program ALM solver ({\bf QPALM})}
		\label{alg:QPALM}
	\begin{algorithmic}[1]
	\pretocmd{\State}{\vspace*{3pt}}{}{}%
	\Require
		\begin{tabular}[t]{@{}l@{}}
			\(x^0\in\R^n\),~
			\(y^0\in\R^m\),~
			\(\varepsilon,\varepsilon_{\rm a},\varepsilon_{\rm r}>0\),~
			\(\vartheta,\rho\in(0,1)\)
		\\
			\(\Delta_\x, \Delta_\y>1\),~
			\(\sigma_\x^{\rm max}, \sigma_\y^{\rm max} > 0\),~
			\(\sigma_\x > 0\)
		\end{tabular}
	\Initialize
		\begin{tabular}[t]{@{}l@{}}
			\(\sigma_\y = \proj_{[10^{-4},10^4]}(20\frac{\max (1, |f(x^0)|)}{\max (1, \|Ax^0 - \proj_C(Ax^0)\|^2)} )\) \cite[\S12.4]{birgin2014practical}%
		\\
			\(\Sx_0=\sigma_\x\I_m\),~
			\(\Sy_0=\sigma_\y\I_m\)%
		\end{tabular}
	\item[{\keywordfont Repeat}]\vspace*{3pt} for \(k=0,1,\ldots\)
	\def\myVar#1{\fillwidthof[l]{y}{#1}}
	\State
		\(x\gets x^k\).
		Let \(\varphi_k\), \(\nabla\varphi_k\) and \(\mathcal J_k\) be as in \eqref{eq:phi}, \eqref{eq:nabla} and \eqref{eq:J}
	\While{ \(\|\nabla\varphi_k(x)\|_{\Sx_k}>\rho^k\varepsilon\) }\label{state:while}%
		\State
			\(\mathtight[0.2]x\gets x+\tau_\star d\) with \(d\) as in \eqref{eq:Newton} and \(\tau_\star\) as in \Cref{alg:LS}%
	\EndWhile
	\State\label{state:xz}%
		\(
			\myVar x^{k+1}=x
		\),~
		\(
			\myVar z^{k+1}=\proj_C(Ax^{k+1}+\Sy_k^{-1}y^k)
		\)
	\State\label{state:y}%
		\(
			\myVar y^{k+1}=y^k+\Sy_k(Ax^{k+1}-z^{k+1})
		\)%
	\State\label{state:if}%
		\If*{$(x^{k+1}, z^{k+1}, y^{k+1})$ satisfies \eqref{subeq:termination}}%
			\State\Return \(x^{k+1}\);%
		\EndIf*%
	\State\label{state:Sigma}%
		\(
			(\Sy_{k+1})_{i,i}
		{}={}
			\begin{cases}[l]
				(\Sy_k)_{i,i} \text{~ if } |(Ax^{k+1}{-{}}z^{k+1})_i|\leq\vartheta|(Ax^k-z^k)_i|,
			\\
				\min\set{\Delta_\y(\Sy_k)_{i,i},\,\sigma_\y^{\rm max}} \text{~ otherwise,}
			\end{cases}
		\)
	\State
		\(
			(\Sx_{k+1})_{i,i} = \min\set{
					\Delta_\x(\Sx_k)_{i,i},\,
					\sigma_\x^{\rm max}
				}
		\)
	\end{algorithmic}
	\end{algorithm}
	
	We now describe the proposed proximal ALM based \Cref{alg:QPALM} for solving \eqref{eq:QP}.
	Steps \ref{state:while}-\ref{state:y} amount to an iteration of proximal ALM.
	As detailed in \Cref{sec:Termination}, the minimization of \(\varphi_k\) needed for the computation of \(x^{k+1}\) can be carried out inexactly.
	In fact, this is done by means of a tailored extremely fast semismooth Newton method with exact line search, discussed in \Cref{sec:Newton,sec:LS}.
	Finally, step \ref{state:Sigma} increases the penalty parameters where the constraint norm has not sufficiently decreased \cite[\S2]{bertsekas1999constrained}.
	
	%
		\subsection{Outer and inner loops: early termination criteria}
			\label{sec:Termination}
	Since the \(x\)-update is not available in closed form, each proximal ALM iteration \((x^k,z^k,y^k)\mapsto(x^{k+1},z^{k+1},y^{k+1})\) in \eqref{eq:PALM} --- which we refer to as an \emph{outer step} --- requires an \emph{inner} procedure to find a minimizer \(x^{k+1}\) of \(\varphi_k\).
	In this subsection we investigate termination criteria both for the outer loop, indicating when to stop the algorithm with a good candidate solution, and for the inner loops so as to ensure that \(x^{k+1}\) is computed with enough accuracy to preserve the convergence guarantees of \Cref{thm:SigmaM}.
	%
	
	\subsubsection{Outer loop termination}
		Following the criterion in \cite{stellato2017osqp}, for fixed absolute and relative tolerances \(\varepsilon_{\rm a},\varepsilon_{\rm r}>0\) we say that \((x,z,y)\) is an \((\varepsilon_{\rm a},\varepsilon_{\rm r})\)-optimal triplet if
		\(
			y
		{}\in{}
			\ncone_C(z)
		\)
		and the following hold:
		\begin{subequations}\label{subeq:termination}
		{\mathtight[0.33]%
		\begin{align} \label{eq:termination1}
			\|Qx+q+\trans Ay\|_\infty
		{}\leq{} &
			\varepsilon_{\rm a}
			{}+{}
			\varepsilon_{\rm r}\max\set{
				\|Qx\|_\infty,~
				\|\trans Ay\|_\infty,~
				\|q\|_\infty
			}
		\\ \label{eq:termination2}
			\|Ax-z\|_\infty
		{}\leq{} &
			\varepsilon_{\rm a}
			{}+{}
			\varepsilon_{\rm r}\max\set{
				\|Ax\|_\infty,~
				\|z\|_\infty
			}.
		\end{align}}%
		\end{subequations}
			From the expression of \(z^{k+1}\) at step \ref{state:xz} it follows that \(Ax^{k+1}+\Sy_k^{-1}y^k-z^{k+1}\in\ncone_C(z^{k+1})\), cf. \eqref{eq:proj}, and hence \(y^{k+1}\) as in step \ref{state:y} satisfies \(y^{k+1}\in\ncone_C(z^{k+1})\).
			A triplet \((x^k,z^k,y^k)\) generated by \Cref{alg:QPALM} is thus \((\varepsilon_{\rm a},\varepsilon_{\rm r})\)-optimal if it satisfies \eqref{subeq:termination}. 
	%
	%


	\subsubsection{Inner loop termination}
		As shown in \Cref{thm:SigmaM}, convergence to a solution can still be ensured when the iterates are computed inexactly, provided that the errors have finite sum.
		Since \(\varphi_k\) as in \eqref{eq:phi} is \(\Sx_k^{-1}\)-strongly convex, from \eqref{eq:PALM} we have that
		\[
			\|\nabla\varphi_k(x)\|_{\Sx_k}
		{}={}
			\|\nabla\varphi_k(x)-\nabla\varphi_k(\tilde x^{k+1})\|_{\Sx_k}
		{}\geq{}
			\|x-\tilde x^{k+1}\|_{\Sx_k^{-1}}.
		\]
		Consequently, condition at \cref{state:while} ensures that \(\|x-\tilde x^{k+1}\|\leq\rho^k\varepsilon'\) holds for all \(k\), with \(\varepsilon'=\varepsilon\sigma_\x^{\rm max}\).
		In turn, \(\|y^{k+1}-\tilde y^{k+1}\|\leq\|A\|\rho^k\varepsilon'\) follows from nonexpansiveness of \(\proj_C\) and finally
		\(
			\|y^{k+1}-\tilde y^{k+1}\|
		{}\leq{}
			2\|\Sy_k\|\|A\|\rho^k\varepsilon'
		{}\leq{}
			2\|\Sy_\infty\|\|A\|\rho^k\varepsilon'
		\).
		As a result, the inner termination criterion at \cref{state:while} guarantees that the error \(\|(x^{k+1},y^{k+1})-(\id+\Sigma_k\M)^{-1}(x^k,y^k)\|\) is summable as required in the convergence analysis of \Cref{thm:SigmaM}.
	In the rest of the section we describe how the \(x\)-update can be carried out with an efficient minimization strategy.
		\subsection{Semismooth Newton method}
			\label{sec:Newton}
	The diagonal matrix \(P_k(x)\) with entries
	\[
		(P_k(x))_{ii}
	{}={}
		\begin{ifcases}
			1 & \ell_i\leq (Ax+\Sy_k^{-1}y^k)_i\leq u_i\\
			0 \otherwise
		\end{ifcases}
	\]
	is an element of the generalized Jacobian \cite[\S7.1]{facchinei2007finite} of \(\proj_C\) at \(Ax+\Sy_k^{-1}y^k\), see \eg \cite[\S6.2.d]{themelis2018acceleration}.
	Consequently, one element \(H_k(x)\in\sym_{++}(\R^m)\) of the generalized Hessian of \(\varphi_k\) at \(x\) is 
	\[
		H_k(x)
	{}={}
		Q+\trans A\Sy_k(\I-P_k(x))A
		{}+{}
		\Sx_k^{-1}.
	\]
	Denoting
	\begin{equation}\label{eq:J}
		\mathcal J_k(x)
	{}\coloneqq{}
		\set{i}[
			(Ax+\Sy_k^{-1}y^k)_i
		{}\notin{}
			{[\ell_i,u_i]}
		],
	\end{equation}
	one has that \((\I-P_k(x))_{ii}\) is 1 if \(i\in\mathcal J_k(x)\) and 0 otherwise.
	Consequently, we may rewrite the generalized Hessian matrix \(H_k(x)\) in a more economic form as
	\begin{equation}\label{eq:H}
		H_k(x)
	{}={}
		Q+\trans A_{\mathcal J_k(x)}(\Sy_k)_{\mathcal J_k(x)}A_{\mathcal J_k(x)}
		{}+{}
		\Sx_k^{-1},
	\end{equation}
	where \(A_{\mathcal J_k(x)}\) is the stacking of the \(j\)-th rows of \(A\) with \(j\in\mathcal J_k(x)\), and similarly \((\Sy_k)_{\mathcal J_k(x)}\) is obtained by removing from \(\Sy_k\) all the \(i\)-th rows and columns with \(i\notin\mathcal J_k(x)\).
	
	A semismooth Newton direction \(d\) at \(x\) solves \(H_k(x)d=-\nabla\varphi_k(x)\).
	Denoting \(\lambda\coloneqq(\Sy_k)_{\mathcal J_k(x)}A_{\mathcal J_k(x)} d\), the computation of \(d\) is equivalent to solving the linear system
	\begin{equation}\label{eq:Newton}
		\begin{bmatrix}
			Q+\Sx_k^{-1} & \trans A_{\mathcal J_k(x)} \\
			A_{\mathcal J_k(x)} & -(\Sy_k)_{\mathcal J_k(x)}^{-1}
		\end{bmatrix}
		\begin{bmatrix}
			d\\
			\lambda
		\end{bmatrix}
	{}={}
		\begin{bmatrix}
			-\nabla\varphi_k(x)\\
			0
		\end{bmatrix}.
	\end{equation}
		\subsection{An exact line search}
			\label{sec:LS}
	\begin{algorithm}[t]
		\algcaption{Exact line search}
		\label{alg:LS}
	\begin{algorithmic}[1]
	\Require
		\(x,d\in\R^n\),~
		diagonal \(\Sigma\in\sym_{++}(\R^n)\)
	\Provide
		optimal stepsize \(\tau_\star\in\R\)
	\State\vspace*{1.55pt}
		Let \(\func{\psi'}{\R}{\R}\), \(\alpha,\beta\in\R\) and \(\delta,\eta\in\R^{2m}\) be as in \eqref{subeq:psi'}
	\State\vspace*{1.55pt}
		Define the set of breakpoints of \(\psi'\)
		
		\noindent\(
			T
		{}={}
			\set{
				\frac{\alpha_i}{\delta_i}
			}[
				i=1,\ldots,2m,~
				\delta_i\neq 0
			]
		\)
	\State\vspace*{1.55pt}
		Sort \(T=\set{t_1,t_2,\ldots}\) such that \(t_i<t_{i+1}\) for all \(i\)
	\State\vspace*{1.55pt}
		Let \(t_i\in T\) be the smallest such that \(\psi'(t_i)\geq0\)
	\State\vspace*{1.55pt}
		\Return
			\(
				\tau_\star
			{}={}
				t_{i-1}
				{}-{}
				\frac{t_i-t_{i-1}}{\psi'(t_i)-\psi'(t_{i-1})}\psi'(t_{i-1})
			\)
			\quad{\footnotesize
				(or \(\tau_\star=t_1\) if \(i=1\))%
			}%
	\end{algorithmic}
	\end{algorithm}
	\begin{subequations}\label{subeq:psi'}%
		\(\nabla\varphi_k(x)\) is piecewise linear, hence so are its sections%
		\[
			\R\ni \tau
		{}\mapsto{}
			\psi_{k,(x,d)}(\tau)
			{}\coloneqq{}
			\varphi_k(x+\tau d)
		\]
		for any \(d\in\R^n\).
		This implies that, given a candidate update direction \(d\), the minimization of \(\varphi_k\) can be carried out \emph{without} the need to perform backtrackings, as an optimal stepsize
		\[
			\tau_\star
		{}\in{}
			\argmin_{\tau\in\R^n}{
				\psi_{k,(x,d)}(\tau)
			}
		\]
		can be explicitly computed, owing to the fact that \(\psi_{k,(x,d)}'\) is a piecewise linear increasing function.
		Indeed, it follows from \eqref{eq:z} and \eqref{eq:nabla} that
		\begin{align*}
			\psi'(\tau)
		{}={} &
			\innprod{\nabla\varphi_k(x+\tau d)}{d}
		\\
		{}={} &
			\innprod{d}{\nabla f(x+\tau d)+\Sx_k^{-1}(x+\tau d -x^k)}
		\\
		&
			{}+{}
			\innprod{Ad}{
				y^k
				{}+{}
				\Sy_k\bigl(
					A(x+\tau d)-Z_k(x+\tau d)
				\bigr)
			}
		\\
		{}={} &
			\mbox{\(\mathtight[0.75]
				\tau\innprod{d}{(Q+\Sx_k^{-1})d}
				{}+{}
				\innprod{d}{Qx+\Sx_k^{-1}(x-x^k)+q}
			\)}
		\\
		&
			\mbox{\(\mathtight[0.75]
				{}+{}
				\innprod{\Sy_k Ad}{
					\bigl[
						Ax+\Sy_k^{-1}y-u
						{}+{}
						\tau Ad
					\bigr]_+
				}
			\)}
		\\
		&
			\mbox{\(\mathtight[0.75]
				{}-{}
				\innprod{\Sy_k Ad}{
					\bigl[
						\ell-Ax-\Sy_k^{-1}y
						{}-{}
						\tau Ad
					\bigr]_+
				}
			\)}
		\\\numberthis\label{eq:psi'}
		{}={} &
			\eta\tau+\beta+\innprod{\delta}{[\delta\tau-\alpha]_+},
		\end{align*}
		where
		\begin{equation}
			\begin{cases}[r @{{}\ni{}} rl]
				\R
				&
				\eta
			{}\coloneqq{}&
				\innprod{d}{(Q+\Sx_k^{-1})d},
			\\
				\R
				&
				\beta
			{}\coloneqq{}&
				\innprod{d}{Qx+\Sx_k^{-1}(x-x^k)+q},
			\\
				\R^{2m}
				&
				\delta
			{}\coloneqq{}&
				\bigl[-\Sy_k^{\nicefrac12}Ad~~~\Sy_k^{\nicefrac12}Ad\bigr],
			\\
				\R^{2m}
				&
				\alpha
			{}\coloneqq{}&
				\Sy_k^{-\nicefrac12}\bigl[y+\Sy_k(Ax-\ell)~~~\Sy_k(u-Ax)\mathrlap{{}-y\bigr].}
			\end{cases}
		\end{equation}
	\end{subequations}
	Due to convexity, it now suffices to find \(\tau\) such that the expression in \eqref{eq:psi'} is zero, as done in \Cref{alg:LS}.
	We remark that,
	since \(\varphi_k\in C^1\) is strongly convex and piecewise quadratic, the proposed nonsmooth Newton method with exact linesearch would converge in \emph{finitely many iterations even with zero tolerance} \cite{sun1997piecewise}.

	\section{Implementation aspects}
		\label{sec:Implementation}
	This section discusses some of the implementation details that are necessary to make QPALM an efficient and competitive algorithm, such as the solution of the linear system at every iteration, preconditioning and infeasibility detection.
	
	\subsection{Linear system}
	We solve the linear system \eqref{eq:Newton} by means of sparse Cholesky factorization routines.
	In the first iteration and after every outer iteration, a sparse LDL factorization of the generalized Hessian matrix \(H_k(x^k)\) as in \eqref{eq:H}
	is computed.
	In between inner iterations, the set of active constraints $\mathcal J_k$ typically does not change much.
	Consequently, instead of doing an LDL factorization from scratch, two low rank updates are sufficient, one for the constraints that enter the active set and one for those that leave.
	As such, the algorithm allows for active set changes where more than one constraint is added and/or dropped, in contrast to active-set methods.
	Therefore, our algorithm typically requires substantially fewer iterations than active-set methods to find the set of constraints that is active at the solution, while still having the advantage of relatively cheap factorization updates.
	The aforementioned routines are carried out with software package {\sc cholmod} \cite{chen2008algorithm}.
	 
	\subsection{Preconditioning} \label{subsec:Preconditioning}
	Preconditioning of the problem data aims at mitigating possible adverse effects of ill conditioning.
	This amounts to scaling problem \eqref{eq:QP} to 
	\begin{equation}\label{eq:scaledQP}
		\minimize_{\bar{x}\in\R^n}\tfrac12\innprod{\bar{x}}{\bar{Q}\bar{x}}+\innprod {\bar{q}}{\bar{x}}
	\quad\stt{}
		\bar{\ell}\leq \bar{A}\bar{x}\leq \bar{u},
	\end{equation}
	with $\bar{x} = D^{-1}x$, $\bar{Q} = c_f DQD$, $\bar{q} = c_f D q$, $\bar{A} = EAD$, $\bar{\ell} = E \ell$ and $\bar{u} =E u$.
	The dual variables in this problem are $\bar{y} = c_f E^{-1} y$.
	Matrices $D\in\R^{n\times n}$ and $E\in\R^{m\times m}$ are diagonal and computed by performing a modified Ruiz equilibration \cite{ruiz2001scaling} on the constraint matrix $A$, scaling its rows and columns so as to have an infinity norm close to 1, as we observed this tends to reduce the number and scope of changes of the active set.
	The scaling factor $c_f$ for scaling the objective was obtained from \cite[\S12.5]{birgin2014practical}, namely
	\(
	c_f = \max (1, \|\nabla f(x_0)\|_{\infty})^{-1}.
	\)
	
	We say that the problem is solved when the unscaled termination criteria \eqref{subeq:termination} holds in a triplet \((\bar{x}^k,\bar{z}^k,\bar{y}^k)\), resulting in the following unscaled criterion
	\begin{align*}
	& c_f^{-1}\|D^{-1}(\bar{Q}\bar{x}^k+\bar{q}+\trans A\bar{y}^k)\|_\infty
				{}\leq{}  \\
					& \hspace{1cm}\varepsilon_{\rm a}
					{}+{}
					\varepsilon_{\rm r}c_f^{-1}\max\set{
						\|D^{-1}\bar{Q}\bar{x}^k\|_\infty,~
						\|D^{-1}\trans{\bar{A}}\bar{y}^k\|_\infty,~
						\|D^{-1}\bar{q}\|_\infty
					},
					\\ 
					&\|E^{-1}(\bar{A}\bar{x}^k-\bar{z}^k)\|_\infty
				{}\leq{}
					\varepsilon_{\rm a}
					{}+{}
					\varepsilon_{\rm r}\max\set{
						\|E^{-1}\bar{A}\bar{x}^k\|_\infty,~
						\|E^{-1}\bar{z}^k\|_\infty
					}.
	\end{align*}

	\subsection{Infeasibility detection}
	The proposed method can also detect whether the problem is primal or dual infeasible from the iterates, making use of the criteria given in \cite{osqp-infeasibility}.
	Let $\delta \bar{y}$ denote the (potential) change in the dual variable, $\delta \bar{y} = \Sy\bigl(\bar{A}\bar{x}-\proj_C(\bar{A}\bar{x}+\Sy^{-1}\bar{y})\bigr)$, then the problem is primal infeasible if for $\delta \bar{y} \neq 0$
	\begin{align*}
		\|D^{-1} \trans{\bar{A}} \delta \bar{y}\|_\infty & {}\leq \varepsilon_{\rm p} \|E\delta \bar{y}\|_\infty,
		\\
		\trans{\bar{u}} [\delta \bar{y}]_+ + \trans{\bar{\ell}} [\delta \bar{y}]_- & {}\leq -\varepsilon_{\rm p} \|E\delta \bar{y}\|_\infty
	\end{align*}
	hold, with $c_f^{-1}E\delta \bar{y}$ the certificate of primal infeasibility.
	Let $\delta \bar{x}$ denote the update in the primal variable, then the problem is dual infeasible if for $\delta \bar{x} \neq 0$ the following conditions hold
	\begin{align*}
	\|D^{-1}\bar{Q}\delta \bar{x}\|_\infty &\leq c_f \varepsilon_{\rm d} \|D\delta \bar{x}\|_\infty, \\
	\trans{\bar{q}} \delta \bar{x} & \leq - c_f \varepsilon_{\rm d} \|D\delta \bar{x}\|_\infty,
	\end{align*}
	\[
		(E^{-1}\bar{A}\delta \bar{x})_i\;
		\begin{ifcases}
			\geq \varepsilon_{\rm d} \|D\delta \bar{x}\|_\infty & \bar{u}_i = +\infty
		\\
			\leq -\varepsilon_{\rm d} \|D\delta \bar{x}\|_\infty & \bar{\ell}_i = -\infty
		\\
			\in [-\varepsilon_{\rm d},\varepsilon_{\rm d}] \|D\delta \bar{x}\|_\infty
			\otherwise.
		\end{ifcases}
	\]
	In that case, $D\delta \bar{x}$ is the certificate of dual infeasibility.
	\section{Numerical simulations}
		\label{sec:Numerical simulations}
	The C implementation of the proposed algorithm was tested for various sets of QPs and benchmarked against state-of-the-art QP solvers: the interior point solver Gurobi \cite{gurobi2018gurobi}, the operator splitting based solver OSQP \cite{stellato2017osqp} and the active-set solver qpOASES \cite{ferreau2014qpoases}.
	The
	first two are also programmed in C, and the third is programmed in C++.
	All simulations were performed on a notebook with Intel(R) Core(TM) i7-7600U CPU @ 2.80GHz x 2 processor and 16 GB of memory.
	The problems are solved to medium accuracy, with the tolerances $\varepsilon_{\rm a}$ and $\varepsilon_{\rm r}$ set to $10^{-6}$ for QPALM and OSQP, as are \texttt{terminationTolerance} for qpOASES and \texttt{OptimalityTol} and \texttt{FeasibilityTol} for Gurobi.
	Apart from the tolerances, all solvers were run with their default options. For QPALM, the default parameters in \Cref{alg:QPALM} are: $\varepsilon = 1$, $\rho = 10^{-1}$, $\theta = 0.25$, $\Delta_\x = \Delta_\y = 10$, $\sigma_\x^{\rm max} = \sigma_\y^{\rm max} = 10^8$, and $\sigma_\y = 10^4$. For the test sets below, the runtime in the figures comprises the average runtime on ten problems of the mentioned problem size or conditioning.
	
	\subsection{Linear programs} \label{subsec:randomLP}
\begin{figure}
		\centering
		\includegraphics[width=7.0cm]{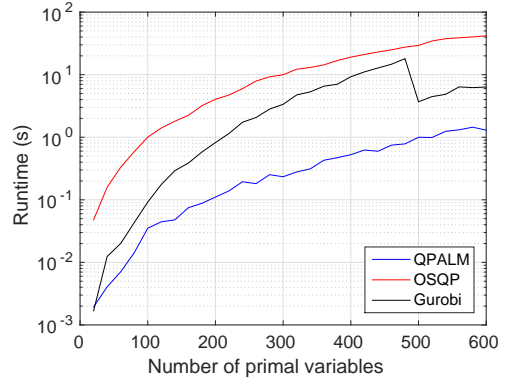} 
		\caption{\small
			Runtime comparison for random LPs of varying sizes.%
		}%
		\label{fig:randomLP}
	\end{figure}

	The first set of tests consists of linear programs (LPs) with randomly generated data.
	Of course, an LP is a special case of a QP with a zero $Q$ matrix.
	The LPs are constructed for 30 values of $n$ equally spaced on a linear scale between 20 and 600. We take $m = 10n$, as typically optimization problems have more constraints than variables.
	The constraint matrix $A \in \R^{m\times n}$ is set to have 50\% nonzero elements drawn from the standard normal distribution, $A_{ij}\sim \mathcal{N}(0,1)$.
	The linear part of the cost $q\in\R^n$ is a dense vector, $q_i \sim \mathcal{N}(0,1)$.
	Finally, the elements of the upper and lower bound vectors, $u,l \in\R^m$, are uniformly distributed on the intervals $[0,1]$ and $[-1,0]$ respectively.
	\Cref{fig:randomLP} illustrates a comparison of the runtimes for QPALM, OSQP and Gurobi applied to random LPs of varying sizes.
	qpOASES is not included in this example, as according to its user manual \cite[\S4.5]{ferreau2014qpoases} it is not suited for LPs of sizes larger than few hundreds of primal variables, which was observed to be the case.
	Also OSQP is not suited for LPs as it hit the maximum number of iterations ($10^5$) for 28 cases and solved inaccurately for the remaining 2 cases.
	QPALM is shown to be an efficient solver for LPs, outperforming the simplex and interior point methods that are concurrently tried by Gurobi.

	\subsection{Quadratic programs} \label{subsec:randomQP}
\begin{figure}
		\centering
		\includegraphics[width=7.0cm]{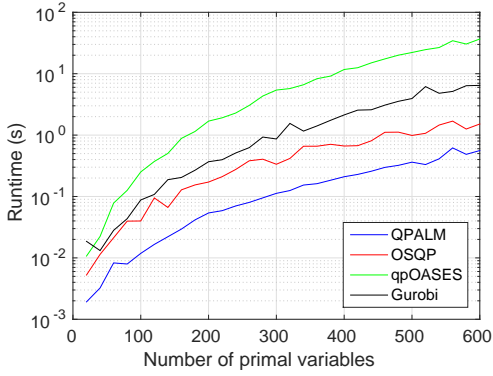} 
		\caption{\small
			Runtime comparison for random QPs of varying sizes.%
		}
		\label{fig:randomQP}
	\end{figure}

	The second set of tests are QPs with data randomly generated in the same way as in \Cref{subsec:randomLP}, with an additional positive definite matrix $Q = M \trans M$, with $M \in \R^{n\times n}$ and 50\% nonzero elements $M_{ij} \sim \mathcal{N}(0,1)$.
	\Cref{fig:randomQP} illustrates the runtimes of the four solvers for such random QPs.
	It is clear that QPALM outperforms the other state-of-the-art solvers regardless of the problem size.

	\subsection{Ill-conditioned problems}
\begin{figure}
		\centering
		\includegraphics[width=7.0cm]{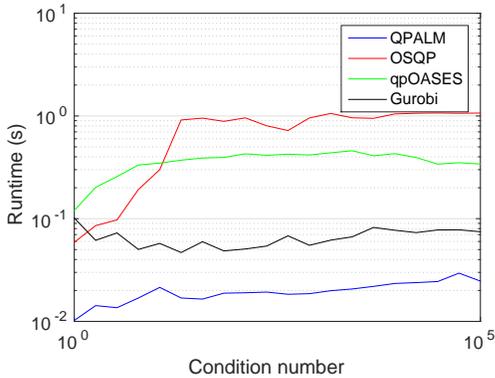} 
		\caption{\small
			Runtime comparison for random QPs with varying conditioning.%
		}%
		\label{fig:randomQPconditioning}
	\end{figure}

	The third test set concerns the conditioning of quadratic programs.
	In this example, the data from \Cref{subsec:randomQP} are reproduced for a QP with $n=100$, but now the impact of the problem conditioning is investigated.
	For this, we set the condition number $\kappa$ of the matrices $Q$ and $A$ (using \texttt{sprandn} in MATLAB), and also scale $q$ with $\kappa$.
	\Cref{fig:randomQPconditioning} shows the runtime results for 20 values of $\kappa$ equally spaced on a logarithmic scale between $10^0$ and $10^5$.
	This figure clearly demonstrates that the first-order method OSQP suffers from ill conditioning in the problem despite the offline Ruiz equilibration it operates.
	From condition number 38 and onwards, OSQP hit the maximum number of iterations ($10^5$).
	Also qpOASES experienced difficulties with ill-conditioned problems.
	From condition number 4833 onwards, it started reporting that the problem was infeasible, while QPALM and Gurobi solved to the same optimal solution.
	From these results it follows that QPALM, supported by preconditioning as discussed in \Cref{subsec:Preconditioning}, is competitive with other solvers in terms of robustness to the scaling of the problem data.
	
	%
	%

	\section{Conclusion}
		\label{sec:Conclusion}
	This paper presented QPALM, a proximal augmented Lagrangian based solver for convex QPs that proved to be efficient and robust against scaling in the problem data.
	Inner minimization procedures rely on semismooth Newton directions and an exact line search which is available in closed form. The iterates, with sparse factorization update routines, allow for large updates in the active set and are more efficient than those of interior point methods and more effective than those of active-set methods. QPALM was shown to compare favorably against state-of-the-art QP solvers, both in runtime and in robustness against problem ill conditioning.
	
	Future work can be focused on considering warm-starting aspects, investigating extensions to nonconvex QPs and SOCPs, and on executing a more thorough set of benchmarking examples focused on problems arising from real applications instead of randomly generated ones.
	

	\bibliographystyle{plain}
	\bibliography{Bibliography_abbr.bib}

\begin{thebibliography}{10}

\bibitem{mosek}
MOSEK ApS.
\newblock {\em Introducing the {MOSEK} {O}ptimization {S}uite 8.1.0.80}, 2019.

\bibitem{axehill2008dual}
D.~Axehill and A.~Hansson.
\newblock A dual gradient projection quadratic programming algorithm tailored
  for model predictive control.
\newblock In {\em 2008 47th IEEE Conference on Decision and Control}, pages
  3057--3064. IEEE, 2008.

\bibitem{osqp-infeasibility}
G.~Banjac, P.~Goulart, B.~Stellato, and S.~Boyd.
\newblock Infeasibility detection in the alternating direction method of
  multipliers for convex optimization.
\newblock {\em optimization-online.org}, 2017.

\bibitem{bauschke2017convex}
H.~H. Bauschke and P.~L. Combettes.
\newblock {\em Convex analysis and monotone operator theory in {H}ilbert
  spaces}.
\newblock CMS Books in Mathematics. Springer, 2017.

\bibitem{bertsekas1999constrained}
D.~P. Bertsekas.
\newblock {\em Constrained optimization and {L}agrange multiplier methods}.
\newblock Athena Scientific, 1999.

\bibitem{birgin2014practical}
E.~G. Birgin and J.~M. Mart\'inez.
\newblock {\em Practical augmented {L}agrangian methods for constrained
  optimization}, volume~10.
\newblock SIAM, 2014.

\bibitem{chen2008algorithm}
Y.~Chen, T.~A. Davis, W.~W. Hager, and S.~Rajamanickam.
\newblock Algorithm 887: {CHOLMOD}, supernodal sparse {C}holesky factorization
  and update/downdate.
\newblock {\em ACM Transactions on Mathematical Software (TOMS)}, 35(3):22,
  2008.

\bibitem{combettes2013variable}
P.~L. Combettes and B.~C. V{\~u}.
\newblock Variable metric quasi-{F}ej{\'{e}}r monotonicity.
\newblock {\em Nonlinear Analysis, Theory, Methods and Applications},
  78(1):17--31, 2013.

\bibitem{Dhingra2017second}
N.~K. Dhingra, S.~Z. Khong, and M.~R. Jovanovi{\'c}.
\newblock A second order primal-dual algorithm for nonsmooth convex composite
  optimization.
\newblock In {\em 2017 IEEE 56th Annual Conference on Decision and Control
  (CDC)}, pages 2868--2873, Dec 2017.

\bibitem{facchinei2007finite}
F.~Facchinei and JS. Pang.
\newblock {\em Finite-dimensional variational inequalities and complementarity
  problems}.
\newblock Springer Science \& Business Media, 2007.

\bibitem{ferreau2008online}
H.~J. Ferreau, H.~G. Bock, and M.~Diehl.
\newblock An online active set strategy to overcome the limitations of explicit
  {MPC}.
\newblock {\em International Journal of Robust and Nonlinear Control},
  18(8):816--830, 2008.

\bibitem{ferreau2014qpoases}
H.~J. Ferreau, C.~Kirches, A.~Potschka, H.~G. Bock, and M.~Diehl.
\newblock {qpOASES}: A parametric active-set algorithm for quadratic
  programming.
\newblock {\em Mathematical Programming Computation}, 6(4):327--363, 2014.

\bibitem{gilbert2014oqla}
J.~C. Gilbert and {\'E}.~Joannopoulos.
\newblock {OQLA}/{QPALM}--{C}onvex quadratic optimization solvers using the
  augmented {L}agrangian approach, with an appropriate behavior on infeasible
  or unbounded problems.
\newblock 2014.

\bibitem{gurobi2018gurobi}
LLC Gurobi~Optimization.
\newblock Gurobi optimizer reference manual, 2018.

\bibitem{nocedal2006numerical}
J.~Nocedal and S.~Wright.
\newblock {\em Numerical optimization}.
\newblock Springer Science \& Business Media, 2006.

\bibitem{parikh2014proximal}
N.~Parikh and S.~Boyd.
\newblock Proximal algorithms.
\newblock {\em Foundations and Trends{\textregistered} in Optimization},
  1(3):127--239, 2014.

\bibitem{patrinos2011global}
P.~Patrinos, P.~Sopasakis, and H.~Sarimveis.
\newblock A global piecewise smooth {N}ewton method for fast large-scale model
  predictive control.
\newblock {\em Automatica}, 47(9):2016--2022, 2011.

\bibitem{rockafellar1976augmented}
R.~T. Rockafellar.
\newblock Augmented {L}agrangians and applications of the proximal point
  algorithm in convex programming.
\newblock {\em Mathematics of operations research}, 1(2):97--116, 1976.

\bibitem{rockafellar2011variational}
R.~T. Rockafellar and R.~JB. Wets.
\newblock {\em Variational analysis}, volume 317.
\newblock Springer Science \& Business Media, 2011.

\bibitem{ruiz2001scaling}
Daniel Ruiz.
\newblock A scaling algorithm to equilibrate both rows and columns norms in
  matrices.
\newblock Technical report, Rutherford Appleton Laboratorie, 2001.

\bibitem{stellato2017osqp}
B.~Stellato, G.~Banjac, P.~Goulart, A.~Bemporad, and S.~Boyd.
\newblock {OSQP}: An operator splitting solver for quadratic programs.
\newblock In {\em 2018 UKACC 12th International Conference on Control
  (CONTROL)}, pages 339--339, Sep. 2018.

\bibitem{sun1997piecewise}
J.~Sun.
\newblock On piecewise quadratic {N}ewton and trust region problems.
\newblock {\em Mathematical Programming}, 76(3):451--467, Mar 1997.

\bibitem{themelis2018acceleration}
A.~Themelis, M.~Ahookhosh, and P.~Patrinos.
\newblock On the acceleration of forward-backward splitting via an inexact
  {N}ewton method.
\newblock In R.~Luke, H.~Bauschke, and R.~Burachik, editors, {\em Splitting
  Algorithms, Modern Operator Theory, and Applications}. Springer.
\newblock To appear \url{https://arxiv.org/abs/1811.02935}.

\end{thebibliography}
\end{document}